\documentclass{elsart}

\usepackage{amssymb,amsbsy}
\usepackage[mathcal]{euscript}

\def\bphi{\bar{\phi}}
\def\bchi{\bar{\chi}}
\def\sphi{\phi^*}
\def\schi{\chi^*}
\def\rmg{{\rm g}} 

\def\nb{\nabla}
\def\bnb{\bar{{\nabla}}}
\def\bR{\bar{R}}
\def\La{\Lambda}
\def\Up{\Upsilon}

\def\bUp{\bar{\Up}}
\def\bLa{\bar{\La}}

\def\z{\theta}

\def\fr{\frac}
\def\d{\partial}

\def\g{\gamma}
\def\Chr{\Gamma}

\def\bbg{\bar{\gamma}}
\def\vecX{{\mathfrak X}}
\def\xiv{\stackrel{\mbox{\tiny $\rightarrow$}}{\boldsymbol\xi}}
\def\lie{\pounds_{\xiv}}

\def\be{\begin{equation}}
\def\ee{\end{equation}}
\def\bnr{\begin{eqnarray*}}
\def\enr{\end{eqnarray*}}
\def\bea{\begin{eqnarray}}
\def\eea{\end{eqnarray}}

\begin{document}

\begin{frontmatter}



\title{Bi-conformal vector fields and the local geometric 
characterization of conformally separable pseudo-Riemannian manifolds II}

\author{Alfonso Garc\'{\i}a-Parrado G\'omez-Lobo}

\address{Departamento de F\'{\i}sica Te\'orica, 
Universidad del Pa\'{\i}s Vasco. Apartado 644, 48080 Bilbao 
(Spain)}

\ead{wtbgagoa@lg.ehu.es}

\begin{abstract}
In this paper we continue the study of bi-conformal vector fields 
started in {\em Class. Quantum Grav.} {\bf 21} 2153-2177.
These are vector fields defined on a pseudo-Rieman\-ni\-an 
manifold by the differential conditions 
$\lie P_{ab}=\phi P_{ab}$, $\lie\Pi_{ab}=\chi\Pi_{ab}$ 
where $P_{ab}$, $\Pi_{ab}$ are
 orthogonal and complementary projectors 
with respect to the metric tensor $\rmg_{ab}$.
In a previous paper we explained how the analysis of these 
differential conditions enabled us to derive local
geometric characterizations of the most relevant cases of 
{\em conformally separable} (also called double 
twisted) pseudo Riemannian manifolds.
In this paper we carry on this analysis further and 
provide local invariant characterizations of 
conformally separable pseudo-Riemannian manifolds 
with {\em conformally flat} leaf metrics. 
These characterizations are rather similar 
to that existing for conformally flat pseudo-Riemannian manifolds 
but instead of the Weyl tensor, we must demand the vanishing of 
certain four rank tensors constructed from the curvature of 
an affine, non-metric, connection (bi-conformal connection). 
We also speculate with possible applications to finding
results for the existence of foliations by conformally flat 
hypersurfaces in any pseudo-Riemannian manifold.
\end{abstract}

\end{frontmatter}

\section{Introduction}
\label{intro}
One of the most interesting problems faced in Differential Geometry is the 
invariant characterization of pseudo-Riemannian manifolds. 
The most famous exposition about this issue was stated
by Klein in his {\em Erlangen program} \cite{ERLANGEN} where he set about  
to the ambitious task of classifying all the possible geometries be they 
underlay by a pseudo-Riemannian manifold or not. 
  
If we stick to pseudo-Riemannian manifolds then the invariant 
characterizations are usually achieved by means of the 
definition of intrinsic (coordinate independent) geometric objects  
upon which a condition is imposed. The best known examples are flat 
and conformally flat 
pseudo-Riemannian manifolds which are those such that the Riemann tensor and 
the Weyl tensor 
are zero respectively. Other relevant cases found in the framework of General Relativity 
are the Schwarzschild geometry \cite{FERRANDO}, 
Kerr black hole \cite{SIMON,MARS},
plane fronted waves with parallel rays 
and many others which are invariantly characterized by certain geometric 
conditions. 

Symmetry considerations often play an important role in invariant 
characterizations (in fact they lie at the heart of the Erlangen program). 
They come in the form of concrete Lie groups acting on the manifold 
under study (group realizations)
and by knowing the orbits or the isotropy subgroups of such actions
it is possible sometimes to identify the pseudo-Riemannian manifold.
       
The above description of symmetries is performed 
in terms of {\em finite groups} but we can also 
study the {\em infinitesimal generators} of these
symmetry groups. These are vector fields satisfying certain differential 
conditions which in general involve the Lie derivative. For instance 
in the case of isometries and conformal transformations the 
corresponding differential conditions fulfilled by the infinitesimal generators
are 
\be
\lie\rmg_{ab}=0,\ \ (\mbox{isometries}),\ \lie\rmg_{ab}=2\phi\rmg_{ab},\ 
(\mbox{conformal motions}),
\label{conformal-motions}
\ee
where $\phi$ is a smooth function.
There are other examples of symmetries studied in the literature 
(see e.g. \cite{LEVINE,HALL,EXACT}) all of them involving the basic objects in 
Differential Geometry (Levi-Civita connection, Riemann curvature, 
Ricci tensor, etc).
However, not many examples have been tackled with more general differential 
conditions (some of them can be found in 
\cite{sergi,KERR-SCHILD,TSAMPARLIS}). In any case from the study of the 
differential conditions satisfied by the vector fields generating the 
symmetry we can sometimes get the geometric conditions fulfilled by 
the pseudo-Riemannian manifolds admitting the symmetry under study 
and these conditions are precisely the geometric
characterizations alluded to above. 

In all the cases we are aware of, the differential conditions are linear in
the generating vector fields and one can show that these vector fields are 
(local) Lie algebras. If such Lie algebras are finite dimensional then it 
is possible to derive from the differential conditions the  
{\em normal form} which generically looks like
$$
D_{x^a}\Phi_B=f(x^1,\dots,x^n,\Phi^1,\dots,\Phi^m), 
$$ 
where $\{x^1,\dots,x^n\}$ are local coordinates, $D_{x^a}$, $a=1,\dots n$ 
a differential 
operator and and $\Phi^1,\dots,\Phi^m$ a set of variables called 
the {\em system variables}.  
Particularly interesting
is the study of the {\em first integrability} and {\em complete integrability} 
conditions associated to the normal form because 
they are the geometric conditions under which 
a pseudo-Riemannian manifold admits locally a finite dimensional
Lie algebra of infinitesimal generators with the highest possible dimension. 
For instance a $n$-dimensional  
locally conformally flat pseudo-Riemannian manifold 
with a $C^2$ conformal factor 
can be characterized by the existence of a local
Lie algebra of conformal motions
of dimension $(n+1)(n+2)/2$ and the complete integrability conditions computed
from (\ref{conformal-motions}) are just the vanishing of 
the Weyl tensor ($n>3$) in a neighbourhood of a point. 
Therefore from the sole symmetry considerations 
one can obtain important geometric objects and this is one of our
aims in this paper where we carry out the calculation of the 
first and the complete integrability conditions for the 
case of {\em bi-conformal vector fields} (see next).

In \cite{BI-CONFORMAL} a new symmetry transformation called 
bi-conformal transformation was put forward. The infinitesimal 
generators of such transformations are called bi-conformal vector fields 
and they are defined through the differential conditions 
$$
\lie P_{ab}=\phi P_{ab},\ \ 
\lie\Pi_{ab}=\chi\Pi_{ab},
$$ 
where $P_{ab}$ and $\Pi_{ab}$ are orthogonal and complementary 
projectors with respect to the metric $\rmg_{ab}$. This 
study was improved in \cite{PAPERI} where we obtained a simple 
expression of the normal form
associated to the above conditions and obtained invariant, 
coordinate free, characterizations of the most relevant 
cases of {\em conformally separable} pseudo-Riemannian manifolds
(we will quote this paper very often in this work and so it 
will be referred to as paper I). A pseudo-Riemannian
manifold is conformally separable at a point $q$ if there exist 
a local coordinate chart $x=\{x^1,\dots,x^n\}$
based at $q$ such that the metric tensor takes the form
\be
\hspace{-1cm}ds^2=\Xi_1(x)G_{\alpha\beta}dx^{\alpha}dx^{\beta}+
\Xi_2(x)G_{AB}dx^Adx^B,\ 1\leq\alpha,\beta\leq p,\ 
1\leq A,B,\leq n-p,
\label{adapted}
\ee
where the functions $G_{\alpha\beta}$, $G_{AB}$ only depend on
the coordinates labelled by their respective subindexes (see 
definition \ref{separable}). The metrics $\Xi_1(x)G_{\alpha\beta}$, 
$\Xi_2(x)G_{AB}$ are called {\em leaf metrics} of the separation.
If $G_{\alpha\beta}$, $G_{AB}$ are flat metrics then  
the manifold is {\em bi-conformally flat} at $q$.

This paper is the continuation of the work started in paper I. There 
we introduced a new affine connection (bi-conformal connection) 
and explained how its use allows us to obtain relevant 
geometric information out of the study bi-conformal vector fields. 
In concrete terms we showed that a pseudo-Riemannian manifold is conformally 
separable at a point if and only if a certain rank three tensor 
$T_{abc}$ constructed from a pair 
of orthogonal and complementary projectors is zero in a neighbourhood
of that point. These projectors
are naturally identified with the leaf metrics (see section \ref{bi-flat}
for more details about this identification). 
In this work a local geometric characterization for bi-conformally flat spaces 
in terms of a certain tensor $T^a_{\ bcd}$  constructed also from 
the projectors $P_{ab}$ and $\Pi_{ab}$ is derived.
We prove the remarkable result that both $T_{abc}$, $T^a_{\ bcd}$ are zero 
if and only if the space is bi-conformally flat at a point being $P_{ab}$ 
and $\Pi_{ab}$ the leaf metrics \footnote{If any of the projectors 
have algebraic rank three then the conditions are slightly different} 
This is the translation to the bi-conformal case of the 
familiar condition that the Weyl tensor vanish  
for a pseudo-Riemannian manifold to be locally
conformally flat. The afore-mentioned conditions will be encountered
here as the complete integrability conditions associated to 
the differential conditions defining a bi-conformal vector field.
We can also give geometric conditions under which a conformally separable 
pseudo-Riemannian 
manifold admits conformally flat leaf metrics. These conditions take the form 
$$
\ P^a_{\ r}P_b^{\ s}P_c^{\ q}P_d^{\ t}T^r_{\ sqt}=0.
$$
One of the most remarkable advantages of these characterizations is 
that they enable us to define intrinsically when a pseudo-Riemannian 
manifold is conformally separable and bi-conformally flat at a point.
Furthermore we can test if two given orthogonal and complementary
projectors $P_{ab}$, $\Pi_{ab}$ give rise to a conformally 
separable pseudo-Riemannian manifold with conformally flat leaf metrics
without actually finding the adapted local coordinates of 
(\ref{adapted}).

The paper outline is as follows: in section \ref{preliminaries} we 
recall basic issues of paper I used in this work 
(the paper is self-contained
and there is no need of reading paper I to understand 
the results presented here). After this the paper is divided into two 
differenced parts: in the first one comprised by sections 
\ref{first-integrability} and \ref{complete} we have put all the 
issues dealing strictly with bi-conformal vector fields
being these the first and the complete integrability conditions 
respectively. In the second part we show how these conditions  
can be used to characterize locally bi-conformally flat pseudo
Riemannian manifolds and conformally separable 
pseudo-Riemannian manifolds with conformally 
flat foliations (section \ref{bi-flat}).
These results are gathered by theorems 
\ref{bi-conformal-char} and \ref{bi-conformal-char-2} (the special 
case of any of the leaf metrics being of rank 3 is treated in
theorem \ref{case-3}).
A reader only interested in these geometric characterizations 
should jump straight to this section and skip the long tensor 
calculations of sections \ref{first-integrability} and 
\ref{complete}. 
Finally examples are provided in section \ref{examples}.  

Part of the results presented in this paper relies on hefty tensor
calculations. These calculations have been done by hand and double-checked 
with the newly released 
Mathematica package ``xtensor" \cite{XTENSOR} with excellent 
agreement.  

\section*{Notation conventions}
The notation of the paper is standard. We work in a $C^{\infty}$ connected
pseudo-Riemannian  manifold $V$ with metric tensor $\rmg_{ab}$ and
we use index notation for all objects constructed 
from the tensor bundles $T^r_s(V)$ of $V$. Square brackets enclosing 
indexes are used to denote antisymmetrization and whenever a set 
of indexes is between strokes it is 
excluded from the antisymmetrization
operation. The metric tensor gives 
rise to the Levi-Civita connection 
$\g^a_{\ bc}$, (we reserve the nomenclature $\Gamma^a_{\ bc}$ for 
the connection components 
calculated in a natural basis) and the curvature tensor $R^a_{\ bcd}$, 
being our convention
for the relation between these two
\be
R^a_{\ bcd}\equiv\d_c\Chr^a_{\ db}-\d_d\Chr^a_{\ cb}+
\Chr^a_{\ rc}\Chr^r_{\ db}-\Chr^a_{\ rd}\Chr^r_{\ cb}.
\label{convention}
\ee
Under this convention the Ricci identity becomes 
$$
\nb_b\nb_cu^a-\nb_c\nb_bu^a=R^a_{\ rbc}u^r,\ \ \nb_b\nb_cu_a-\nb_c\nb_bu_a=-R^r_{\ abc}u_r,
$$ 
where $\nb_a$ is the covariant derivative of the Levi-Civita connection.

All the above relations are still valid if $\g^a_{\ bc}$ is a connection with 
no torsion (symmetric or affine connection). 

The infinite dimensional Lie algebra of smooth vector fields of the manifold
$V$ is denoted by $\vecX(V)$. 
Finally the Lie derivative operator 
with respect to a vector field $\xiv$ is $\lie$.

\section{Preliminaries}
\label{preliminaries}
In this section we review concepts of paper I which are needed in this work. 
We quote the results without proofs as they can all be found in paper I.
\begin{defn}
A smooth vector field $\xiv$ on $V$ is said to be a 
{\em bi-conformal vector field} if it fulfills the condition
\be 
\lie P_{ab}=\phi P_{ab},\ \ \ 
\lie\Pi_{ab}=\chi\Pi_{ab},\ \phi,\ \chi\in C^{\infty}(V), 
\label{bi-conformal}
\ee
where $P_{ab}$, $\Pi_{ab}$ are smooth sections of the tensor bundle $T^0_2(V)$ satisfying the properties
\bea
P_{ab}=P_{ba},\ \Pi_{ab}=\Pi_{ba},\ P_{ab}+\Pi_{ab}=\rmg_{ab},\nonumber\\ 
P_{ap}P^{p}_{\ b}=P_{ab},\ \Pi_{ap}\Pi^{p}_{\ b}=\Pi_{ab},\  
P_{ap}\Pi^p_{\ b}=0.
\label{properties}
\eea
\label{Bi-conformal}
\end{defn}

The geometrical meaning of these conditions is that $P_{ab}$ and 
$\Pi_{ab}$ are {\em orthogonal projectors} with respect to the 
metric tensor $\rmg_{ab}$ at each point of the manifold. 
We describe next briefly some useful properties of $P_{ab}$ and
$\Pi_{ab}$ which will be used along the paper.
To start with note that 
we can decompose the vector space $T_p(V)$ as a direct sum 
of the ranges of the endomorphisms $P^a_{\ b}$ and $\Pi^a_{\ b}$. 
These ranges are the respective eigenspaces of the endomorphisms 
with eigenvalue $+1$ and they are orthogonal subspaces. 
$P_{ab}$ and $\Pi_{ab}$  enable us to represent 
any {\em nondegenerate} smooth distribution $D$ in a compact way 
\footnote{By nondegenerate we mean 
that the scalar product $\rmg_{ab}$ restricted to the subspace of $T_p(V)$ 
generated by the distribution $D$ is not degenerated.}. To see this more 
clearly, let $\{u^a_1,\dots,u^a_p\}$, $0<p<n$, be an orthonormal set 
of smooth vector fields spanning the distribution $D$ 
(such a set always exists if $D$ is nondegenerate) with 
$\epsilon_{\alpha}=u^a_{\alpha}\rmg_{ab}u^b_{\alpha}$.  
Then the smooth sections
\be
P_{ab}\equiv\sum_{\alpha=1}^p\epsilon_\alpha u^\alpha_au^\alpha_b,\ 
\Pi_{ab}\equiv\rmg_{ab}-P_{ab},
\label{distribution}
\ee 
satisfy (\ref{properties}). Conversely, 
the smoothness of $P^a_{\ b}$ and $\Pi^a_{\ b}$ together with 
(\ref{properties}) gua\-ran\-tee that the ranges of $P^a_{\ b}$ and 
$\Pi^a_{\ b}$ span smooth distributions on $V$. 
To see this we need to show that 
the range dimension of 
each projector does not vary in the manifold $V$ 
(these numbers are $p\equiv P^a_{\ a}$ and $n-p\equiv\Pi^a_{\ a}$) 
because in that case 
such ranges are smooth distributions (see e.g. \cite{HALL-RENDALL}). 
Let us split \footnote{We are indebted to Graham Hall for this proof.} 
$V$ in subsets $A_k$, $B_{k'}$ $k,k'=1,\dots,n-1$ defined by the conditions
$$
A_{k}=\{q\in V:\ P^a_{\ a}=k\},
\ B_{k'}=\{q\in V:\ \Pi^a_{\ a}=k'\},
$$
and denote by $k_{min}$ and $k_{max}$ respectively 
the minimum and maximum value of the integer $k$ (similarly 
we define  $k'_{min}$, $k'_{max}$). 
By the rank theorem (see e. g. theorem 3.1 of \cite{HALL}) 
$A_{k_{max}}$ and $B_{k'_{max}}$ are 
open sets and from the third property of (\ref{properties}) 
clearly $A_{k_{max}}=B_{k'_{min}}$, $B_{k'_{max}}=A_{k_{min}}$. 
On the other hand using again the rank theorem we 
deduce that $V\setminus B_{k'_{min}}$, $V\setminus A_{k_{min}}$ 
are open and thus 
$B_{k'_{min}}$ and $A_{k_{min}}$ must be closed from which 
we conclude that $A_{k_{max}}$, $B_{k'_{max}}$ are both open and closed 
at the same time and thus equal to $V$ since it is connected.

The differential conditions (\ref{bi-conformal}) 
are the starting point for an interesting study of the properties 
and geometric significance of bi-conformal vector fields. In paper 
I we argued that the set of bi-conformal vector fields of a 
pseudo-Riemannian manifold $V$ is a Lie subalgebra of $\vecX(V)$ and we 
established the conditions under which such algebra is always 
finite dimensional as well
as its greatest dimension $N$. The Lie algebra is finite dimensional 
if $p,\ n-p\neq 1,2$ and in this case 
$$
N=\fr{1}{2}(p+1)(p+2)+\fr{1}{2}(n-p+1)(n-p+2).
$$
The number $N$ is calculated from the {\em normal form} associated with the differential conditions
(\ref{bi-conformal}). This is a set of equations obtained from the differential conditions by means
of successive differentiations and they were calculated in paper I 
(here, this normal form is recalled in (\ref{normal-form})). 
In doing this calculation the introduction of 
a new affine connection ({\em bi-conformal connection}) revealed itself essential rendering the normal form
very neatly. The components of the bi-conformal connection $\bbg^a_{\ bc}$ are related to the 
Levi-Civita connection $\g^a_{\ bc}$ by the relation
$$
\bbg^a_{\ bc}=\g^a_{\ bc}+\fr{1}{2p}(E_bP^a_{\ c}+E_cP^a_{\ b})+\fr{1}{2(n-p)}(W_b\Pi^a_{\ c}+W_c\Pi^a_{\ b})+
\fr{1}{2}(P^a_{\ p}-\Pi^a_{\ p})M^p_{\ bc},   
$$
where
\be
M_{abc}\equiv\nb_{b}P_{ac}+\nb_{c}P_{ab}-\nb_aP_{bc},\ E_a\equiv M_{abc}P^{bc},\ W_a\equiv-M_{abc}\Pi^{bc}.
\label{m-definition}
\ee
By definition the bi-conformal connection is an affine connection and we shall denote the 
covariant derivative 
and the curvature tensor of this connection by $\bnb$ and $\bR^a_{\ bcd}$ respectively. The bi-conformal 
connection 
does not stem from a metric tensor in general as will be shown in explicit examples
and hence the tensor $\bR^a_{\ bcd}$ does not fulfill the same properties 
as the curvature tensor of a metric connection.
The Bianchi identities though, remain as in the case of a metric connection.

In paper I, we explained the role of the bi-conformal connection in the local geometric 
characterization of 
{\em conformally separable} pseudo-Riemannian manifolds (see definition \ref{separable}).
This role will be strengthened in section \ref{bi-flat} 
where we will develop a local invariant characterization of conformally separable 
pseudo-Riemannian manifolds with conformally flat leaf metrics ({\em bi-conformally flat} 
pseudo-Riemannian manifolds). 

\section{First integrability conditions}
\label{first-integrability}
As we commented before, one can differentiate equation (\ref{bi-conformal})
a certain number of times and then isolate the derivatives of certain variables
(system variables) in terms of themselves thereby obtaining a ``closed'' or  
normal form. This calculation was accomplished in paper I and 
we reproduce next the result.
\bea
(a)\ \ \bnb_a\phi&=&\bphi_a+\sphi_a,\ \bnb_a\chi=\bchi_a+\schi_a,\nonumber\\
(b)\ \ \bnb_b\sphi_a&=&\fr{-1}{p}\left[\lie(\bnb_bE_a)+\fr{1}{2}(\bchi_bE_a+\bchi_aE_b-(\bchi^rE_r)\Pi_{ab})\right],\nonumber\\
(c)\ \ \bnb_b\schi_a&=&\fr{1}{p-n}\left[\lie(\bnb_bW_a)+\fr{1}{2}(\bphi_bW_a+\bphi_aW_b-(\bphi^rW_r)P_{ab})\right],\nonumber\\
(d)\ \ \bnb_b\bphi_c&=&\fr{1}{2-p}\left[\lie L^0_{bc}+2\bphi^r\bnb_rP_{bc}\right],\label{normal-form}\\
(e)\ \ \bnb_b\bchi_c&=&\fr{1}{2-n+p}\left[\lie L^1_{bc}+2\bchi^r\bnb_r\Pi_{bc}\right],\nonumber\\
(f)\ \ \bnb_b\xi^a&=&\Psi_b^{\ a},\nonumber\\ 
(g)\ \ \bnb_b\Psi_c^{\ a}&=&\fr{1}{2}(\bphi_bP^a_{\ c}+\bphi_cP^a_{\ b}-\bphi^aP_{cb}+
\bchi_b\Pi^a_{\ c}+\bchi_c\Pi^a_{\ b}-\bchi^a\Pi_{cb})-\xi^d\bar{R}^a_{\ cdb},\nonumber 
\eea
with the definitions 
$$
\sphi_a\equiv \Pi_{ab}\phi^b,\ \bphi_a\equiv P_{ab}\phi^b,\ 
\schi_a\equiv P_{ab}\chi^b,\ \bchi_a\equiv\Pi_{ab}\chi^b,
$$
\bea
\hspace{-1cm}L^0_{\ bc}&&\equiv 
2\left[P^d_{\ r}\bR^r_{\ cdb}-\fr{1}{p}(P^d_{\ c}P^r_{\ q}\bR^q_{\ rdb}+P^d_{\ b}P^r_{\ q}\bR^q_{\ rdc}-
P^r_{\ q}\bR^q_{\ rbc})\right]+\fr{\bR^0}{1-p}P_{bc},\nonumber\\
\bR^0&&\equiv P^d_{\ r}\bR^r_{\ cdb}P^{cb},
\label{L0}
\eea 
\bea
L^1_{\ bc}\equiv 
2\left[\Pi^d_{\ r}\bR^r_{\ cdb}-\fr{1}{n-p}(\Pi^d_{\ c}\Pi^r_{\ q}\bR^q_{\ rdb}+\Pi^d_{\ b}\Pi^r_{\ q}\bR^q_{\ rdc}-
\Pi^r_{\ q}\bR^q_{\ rbc})\right]&+&\nonumber\\
+\fr{\bR^1}{1-n+p}\Pi_{bc},\ \ \bR^1\equiv\Pi^d_{\ r}\bR^r_{\ cdb}\Pi^{cb}.& &
\label{L1}
\eea

The variables lying on the l.h.s. of (\ref{normal-form}) 
are the system variables. 
Not all these variables are independent because as shown in paper I
they are constrained by the conditions (constraint equations) 
\be
(A)\left.\begin{array}{c}
\lie P_{ab}=\phi P_{ab}\\
\lie\Pi_{ab}=\chi\Pi_{ab}
\end{array}\right\},\ \ 
(B)\left.\begin{array}{l}
\lie E_a=-p\sphi_a\\
\lie W_a=-(n-p)\schi_a\end{array}\right\}.
\label{ligaduras}
\ee 
In this section we spell out the first integrability conditions 
of each equation of (\ref{normal-form}). 
These are geometric 
conditions arising from the compatibility conditions yielded by 
the commutation 
of two covariant derivatives. The commutation rules for these 
derivatives are given by the Ricci identity
$$
\hspace{-.3cm}\bnb_a\bnb_b\Xi^{a_1\dots a_r}_{\ b_1\dots b_s}-
\bnb_b\bnb_a\Xi^{a_1\dots a_r}_{\ b_1\dots b_s}
=\sum_{q=1}^r\bR^{a_q}_{\ tab}\Xi^{a_1\dots a_{q-1}ta_{q+1}\dots a_r}_{\ b_1\dots b_s}-
\sum_{q=1}^{s}\bR^t_{\ b_qab}\Xi^{a_1\dots a_r}_{\ b_1\dots b_{q-1}tb_{q+1}\dots b_s},
$$
where we must replace the tensor $\Xi^{a_1\dots a_r}_{\ b_1\dots b_s}$ by the system 
variables and apply 
(\ref{normal-form}) to work out the covariant derivatives. 
Each one of the equations derived in this fashion only 
involves system variables and it
is called {\em first integrability condition}. 
These integrability conditions can be further covariantly differentiated 
yielding integrability conditions of higher degree. The constraint equations 
(\ref{ligaduras}) also give rise to integrability conditions when 
differentiated in the obvious way.

\subsection*{Equation (\ref{normal-form})-f}
This is the simplest integrability condition being its expression 
$$
\bnb_c\bnb_b\xi^a-\bnb_b\bnb_c\xi^a=\bR^a_{\ rcb}\xi^r=\bnb_c\Psi_b^{\ a}
-\bnb_b\Psi_c^{\ a},
$$
which is an identity as is easily checked by replacing the covariant derivatives 
of $\Psi_a^{\ b}$. 
\subsection*{Equation (\ref{normal-form})-g}
The integrability conditions of (\ref{normal-form})-{\em g} are given by 
\bnr
\bnb_a\bnb_b\Psi_c^{\ d}-\bnb_b\bnb_a\Psi_c^{\ d}=
-\Psi_a^{\ r}\bR^d_{\ crb}+\Psi_b^{\ r}\bR^d_{\ cra}-\xi^r\bnb_r\bR^d_{\ cab}+\\
+\fr{1}{2}\nb_{a}((\bphi_bP^d_{\ c}+\bphi_cP^d_{\ b}-\bphi^dP_{cb}+
\bchi_b\Pi^d_{\ c}+\bchi_c\Pi^d_{\ b}-\bchi^d\Pi_{cb}))-\\
-\fr{1}{2}\nb_{b}((\bphi_aP^d_{\ c}+\bphi_cP^d_{\ a}-\bphi^dP_{ca}+
\bchi_a\Pi^d_{\ c}+\bchi_c\Pi^d_{\ a}-\bchi^d\Pi_{ca})).
\enr
We apply now the Ricci identity to the left hand side of this expression and gather 
all the terms containing contractions with the tensor $\Psi_a^{\ b}$ in a single term
by means of the identity
$$
\lie\bR^d_{\ cab}=\xi^r\bnb_r\bR^d_{\ cab}-\Psi_r^{\ d}\bR^r_{\ cab}+
\Psi_c^{\ r}\bR^d_{\ rab}+\Psi_{a}^{\ r}\bR^d_{\ crb}+
\Psi_b^{\ r}\bR^d_{\ car},
$$
from which we obtain 
\bea
\lie\bR^d_{\ cab}=\bnb_{[a}\bphi_{b]}P^d_{\ c}&+&P^d_{\ [b}\bnb_{a]}\bphi_c-P_{c[b}\bnb_{a]}\bphi^d+
\bnb_{[a}\bchi_{b]}\Pi^d_{\ c}+\Pi^d_{\ [b}\bnb_{a]}\bchi_c-\nonumber\\
-\Pi_{c[b}\bnb_{a]}\bchi^d+\bphi_{[b}\bnb_{a]}P^d_{\ c}&+&\bphi_c\bnb_{[a}P^d_{\ b]}-\bphi^d\bnb_{[a}P_{b]c}+
\bchi_{[b}\bnb_{a]}\Pi^d_{\ c}+\bchi_c\bnb_{[a}\Pi^d_{\ b]}-\nonumber\\
-\bchi^d\bnb_{[a}\Pi_{b]c}\label{lie-curvatura}.
\eea
The direct substitution of $\bnb_a\bchi_b$, $\bnb_a\bphi_b$ by the 
expressions given by
(\ref{normal-form})-{\em d} and (\ref{normal-form})-{\em e} respectively  
yields after lengthy algebra
\bea
\fr{1}{2}\lie T^d_{\ cab}&=&\fr{\bphi_r}{2-p}(P^d_{\ [b}\La^r_{\ a]c}+P^d_{\ q}\Up^{rq}_{\ [b}P_{a]c})+\nonumber\\
+\fr{\bchi_r}{2-(n-p)}
(\Pi^d_{\ [b}\bLa^r_{\ a]c}&+&\Pi^d_{\ q}\bUp^{rq}_{\ [b}P_{a]c})
+\bphi_{[b}\bnb_{a]}P^d_{\ c}+\bphi_c\bnb_{[a}P^d_{\ b]}+\bphi^d\bnb_{[b}P_{a]c}+\nonumber\\
&+&\bchi^d\bnb_{[b}\Pi_{a]c}+
\bchi_{[b}\bnb_{a]}\Pi^d_{\ c}+\bchi_c\bnb_{[a}\Pi^d_{\ b]},
\label{cond-g}
\eea   
where by convenience we introduce the tensors
\bea
\hspace{-1cm}\La^d_{\ bc}\equiv2P^{dr}\bnb_rP_{bc},\ \bLa^d_{\ bc}&=&2\Pi^{dr}\bnb_r\Pi_{bc},
\Up^{sc}_{\ b}\equiv2P^{sr}P^{cq}\bnb_rP_{qb}+(2-p)\bnb_bP^{sc},\ \ \ \ \ \ 
\label{lambda}\\
\bUp^{sc}_{\ b}\equiv2\Pi^{sr}\Pi^{cq}\bnb_r\Pi_{qb}&+&(2-n+p)\bnb_b\Pi^{sc},
\label{upsilon}
\eea 
and
\bea
T^d_{\ cab}\equiv 2\bR^d_{\ cab}-\fr{2}{2-p}(P^d_{\ c}L^0_{[ab]}+P^d_{\ [b}L^0_{a]c}+P_{c[a}L^0_{b]q}P^{qd})+
\nonumber\\
-\fr{2}{2-n+p}(\Pi^d_{\ c}L^1_{[ab]}+\Pi^d_{\ [b}L^1_{a]c}+\Pi_{c[a}L^1_{b]q}\Pi^{qd}).
\label{tensor-tt}
\eea
This last tensor will play an important role in the local characterization of 
bi-conformally flat pseudo-Riemannian manifolds as will be seen later.

An interesting invariance property of some of the above tensors needed in future calculations is
\be
\lie\La^d_{\ bc}=0,\ \ \lie\bLa^d_{\ bc}=0
\label{inv-3}
\ee
which are easily obtained from (\ref{cond-ligaduras-up})
(see below).

\subsection*{(\ref{normal-form})-b and (\ref{normal-form})-c} 
Again the calculations are tedious but straightforward. The covariant derivatives of $\bphi_a$ and $\bchi_a$ 
are calculated through (\ref{normal-form})-{\em d} and (\ref{normal-form})-{\em e}
and  to commute the Lie derivative and the 
covariant derivative when differentiating both equations we use the identity (see \cite{YANO,SCHOUTEN})
\bea
\bnb_c\lie T^{a_1\dots a_s}_{\ b_1\dots b_q}
-\lie\bnb_cT^{a_1\dots a_s}_{\ b_1\dots b_q}
=-\sum_{j=1}^s(\lie&\hspace{-1.5cm}\bar{\g}^{a_j}_{cr})T^{\dots a_{j-1}ra_{j+1}\dots}_{\ b_1\dots b_q}+\nonumber\\
&+\sum_{j=1}^q(\lie\bar{\g}^r_{cb_j})T^{a_1\dots a_s}_{\dots b_{j-1}rb_{j+1}\dots},
\label{lie-conmmutation} 
\eea
where, as calculated in paper I, the Lie derivative of the bi-conformal connection is 
$$
\lie\bbg^a_{\ bc}=\fr{1}{2}(\bphi_bP^a_{\ c}+\bphi_cP^a_{\ b}-\bphi^aP_{cb}+\bchi_b\Pi^a_{\ c}+\bchi_c\Pi^a_{\ b}-\bchi^a\Pi_{cb}).
$$
Recall that the Lie derivative of a connection is always a tensor 
even though the connection itself is not (see e.g. \cite{YANO}). 
Putting all this together we get
\bea
& &E_d\lie(\Pi^d_{\ r}T^r_{\ acb})=\bchi_a\bnb_{[c}E_{b]}+\bchi_{[b}\bnb_{c]}E_a-\bchi^r\bnb_{[c}(\Pi_{b]a}E_r)+\nonumber\\
&+&(P^r_{\ a}\bphi_{[c}+\bphi_aP^r_{\ [c}-\bphi^rP_{a[c}+
\Pi^r_{\ a}\bchi_{[c}+\bchi_a\Pi^r_{\ [c}-\bchi^r\Pi_{a[c})\bnb_{b]}E_r+\nonumber\\
&+&\fr{1}{2-n+p}\bchi_qE_r\bUp^{qr}_{\ [b}\Pi_{c]a},\label{cond-b}\\
& &W_d\lie(P^d_{\ r}T^r_{\ acb})=\bphi_a\bnb_{[c}W_{b]}+\bphi_{[b}\bnb_{c]}W_a-\bphi^r\bnb_{[c}(P_{b]a}W_r)+\nonumber\\
&+&(\Pi^r_{\ a}\bchi_{[c}+\bchi_a\Pi^r_{\ [c}-\bchi^r\Pi_{a[c}+
P^r_{\ a}\bphi_{[c}+\bphi_aP^r_{\ [c}-\bphi^rP_{a[c})\bnb_{b]}W_r+\nonumber\\
&+&\fr{1}{2-p}\bphi_qW_r\Up^{qr}_{\ [b}P_{c]a}.\label{cond-c}
\eea
\subsection*{Equation (\ref{normal-form})-a}
The integrability conditions of this equation are easier to handle
\bnr
0=\bnb_a\bnb_b\phi-\bnb_b\bnb_a\phi=\bnb_a\bphi_b-\bnb_b\bphi_a+\bnb_a\sphi_b-\bnb_b\sphi_a\\
0=\bnb_a\bnb_b\chi-\bnb_b\bnb_a\chi=\bnb_a\bchi_b-\bnb_b\bchi_a+\bnb_a\schi_b-\bnb_b\schi_a,\\
\enr
from which we readily obtain by means of (\ref{normal-form})-{\em b}, (\ref{normal-form})-{\em c}, (\ref{normal-form})-{\em d},
(\ref{normal-form})-{\em e}
\bea
\lie\left(\fr{1}{2-p}L^0_{[ab]}-\fr{1}{p}\bnb_{[a}E_{b]}\right)=0,\label{cond-a1}\\
\lie\left(\fr{1}{2-n+p}L^1_{[ab]}-\fr{1}{n-p}\bnb_{[a}W_{b]}\right)=0.\label{cond-a2}
\eea
\subsection*{Equations (\ref{normal-form})-d and (\ref{normal-form})-e}
The integrability conditions yield  
(equation (\ref{inv-3}) is used along the way)
\bea
\fr{2-p}{2}\bphi_rT^r_{\ ceb}&=&2\lie\left(\bnb_{[e}L^0_{b]c}+\fr{1}{2-p}\La^d_{\ c[b}L^0_{e]d}\right)+
\bchi_{[e}L^0_{b]q}\Pi^q_{\ c}+\nonumber\\
+\bchi_c\Pi^q_{\ [e}L^0_{b]q}-\bchi^q\Pi_{c[e}L^0_{b]q}
&+&\fr{2}{2-p}\bphi_r(\La^d_{c[b}\La^r_{e]d}+(2-p)\bnb_{[e}\La^r_{\ b]c}),\label{cond-d}\\
\fr{2-n+p}{2}\bchi_rT^r_{\ ceb}&=&2\lie\left(\bnb_{[e}L^1_{b]c}+\fr{1}{2-n-p}\bLa^d_{\ c[b}L^1_{e]d}\right)+
\bchi_{[e}L^1_{b]q}P^q_{\ c}+\nonumber\\
+\bphi_cP^q_{\ [e}L^1_{b]q}-\bphi^qP_{c[e}L^1_{b]q}&+&
\fr{2}{2-n+p}\bchi_r(\bLa^d_{c[b}\bLa^r_{e]d}+(2-n+p)\bnb_{[e}\bLa^r_{\ b]c}).\label{cond-e} 
\eea
\subsection*{Constraint equations}
Finally only the first integrability conditions coming up from the
constraints (\ref{ligaduras}) are left. These 
integrability conditions 
result from their covariant derivative (we only need to take 
care of (\ref{ligaduras})-{\em A} because the differentiation of (\ref{ligaduras})-{\em B} results in (\ref{normal-form})-{\em b} and 
(\ref{normal-form})-{\em c} which are part of the normal form). 
Differentiation of these equations with respect to $\bnb$
yields after some algebra (apply identity (\ref{lie-conmmutation}))
\be
\lie\bnb_cP_{ab}=\phi\bnb_cP_{ab}+\sphi_cP_{ab},\ \ 
\lie\bnb_c\Pi_{ab}=\chi\bnb_c\Pi_{ab}+\schi_c\Pi_{ab}.
\label{cond-ligaduras-down}
\ee
For completeness we provide also these equations with the $ab$ indexes raised
\be
\lie\bnb_cP^{ab}=-\sphi_cP^{ab}-\phi\bnb_cP^{ab},\ \ 
\lie\bnb_c\Pi^{ab}=-\schi_c\Pi^{ab}-\chi\bnb_c\Pi^{ab},\ \ 
\label{cond-ligaduras-up}
\ee
and the invariance laws
\be
\lie\bnb_cP^a_{\ b}=0,\ \ \lie\bnb_c\Pi^a_{\ b}=0.
\label{mixed-invariance}
\ee
These equations close the whole suite of first integrability conditions. 
In the next section we will obtain geometric information from these
conditions. 

\section{Complete integrability}
\label{complete}
Our next task is to find out when
the first integrability conditions presented in the previous sec\-ti\-on become a set of identities
for every choice of the independent variables of the normal form (\ref{normal-form}). 
This happens if certain geometric conditions (complete integrability conditions) are met. 
Under such conditions there exists in the neighbourhood of each point a finite dimensional 
Lie algebra of bi-conformal vector fields attaining the greatest dimension $N$. 
Thus if we are
able to find pseudo-Riemannian manifolds in which these conditions are 
satisfied we will have proven that the bound $N$ is reached at least locally. 
(see \cite{EISENHARTII}).
The first problem which 
we come across to is that not all the variables appearing in (\ref{normal-form}) 
are independent as there are constraints 
which must be taken into account. However, some of the variables of (\ref{normal-form}) 
are not involved in the constraint
equations (\ref{ligaduras}) and this fact will allow us to find necessary and 
sufficient geometric conditions for the 
whole set of first integrability conditions to become identities. 
The variables which are not constrained by 
(\ref{ligaduras}) are $\bphi_a$ and $\bchi_a$ so we will separate out 
in each of the integrability conditions obtained above all the contributions 
involving these variables. 
If we demand then that $\bphi_a$ and $\bchi_a$ be arbitrary functions 
in the neighbourhood of a point we will obtain a number 
of local geometric conditions 
which will lead us to the complete integrability conditions. 

We will perform next this procedure step by step  
(actually we will not need to analyse all the conditions as some of them 
will turn into identities if the geometric conditions entailed by 
others are imposed). The full calculation is rather cumbersome
and only its main excerpts will be shown so the reader interested 
in the complete integrability conditions should jump directly to theorem 
\ref{complete-integrability}. 

\medskip
\noindent
{\em Equation (\ref{cond-g})}. This equation can be rewritten as 
\be
\lie T^d_{\ cab}=\bphi_rM^{rd}_{\ cab}+\bchi_rN^{rd}_{\ cab}, 
\label{cond-g1}
\ee  
where
\bea
M^{rd}_{\ cab}&\equiv&\fr{2}{2-p}P^r_{\ s}(\La^s_{\ c[a}P^d_{\ b]}+\Up^{sd}_{\ [b}P_{a]c})+2P^r_{\ [b}\bnb_{a]}P^d_{\ c}+\nonumber\\
+2P^r_{\ c}\bnb_{[a}P^d_{\ b]}&+&2P^{dr}\bnb_{[b}P_{a]c},
\label{def-M}\\
N^{rd}_{\ cab}&\equiv&\fr{2}{2-n+p}\Pi^r_{\ s}(\bLa^s_{\ c[a}\Pi^d_{\ b]}+\bUp^{sd}_{\ [b}\Pi_{a]c})+2\Pi^r_{\ [b}\bnb_{a]}\Pi^d_{\ c}+\nonumber\\
+2\Pi^r_{\ c}\bnb_{[a}\Pi^d_{\ b]}
&+&2\Pi^{dr}\bnb_{[b}\Pi_{a]c}.\label{def-N}
\eea
Under the assumption of complete integrability
(\ref{cond-g1}) must be true for every $\bphi_a$, $\bchi_a$ so we have
\be
P^r_{\ s}M^{sd}_{\ cab}=M^{rd}_{\ cab}=0,\ \ \Pi^r_{\ s}N^{sd}_{\ cab}=N^{rd}_{\ cab}=0.
\label{cuarenta}
\ee
Let us study these geometric conditions (it is enough to concentrate on the first condition because the 
other is dual and it is obtained by the usual replacements). Contracting the indexes $d$ and $b$ in this equation we get
\be
\fr{2}{2-p}[pP^{rq}\bnb_qP_{ac}-P^{rq}(P_a^{\ s}\bnb_qP_{sc}+P_c^{\ s}\bnb_qP_{sa})]+P^{sr}\bnb_sP_{ac}=0\label{d-b}.
\ee
Multiplying this by $P_z^{\ a}$ gives the condition
\be
P_z^{\ s}P^{rq}\bnb_qP_{sc}=0,
\ee
which put back in (\ref{d-b}) yields
\be
P^{rs}\bnb_{s}P_{ac}=0\Rightarrow\La^r_{\ ac}=0,\ \ \Up^{sc}_{\ b}=(2-p)\bnb_bP^{sc}.
\label{1st-condition}
\ee
Now using this result we contract $r$ and $b$ in the first of (\ref{cuarenta})
\be
(p-1)\bnb_aP^d_{\ c}-P_a^{\ s}\bnb_sP^d_{\ c}-P_c^{\ s}\bnb_sP^d_{\ a}+P_c^{\ s}\bnb_aP^d_{\ s}=0,
\label{r-b}
\ee
from which we get
\be
P_a^{\ q}P_c^{\ s}\bnb_qP^d_{\ s}=0.
\ee
Combining this with (\ref{r-b}) we readily obtain
\be
\bnb_aP^d_{\ c}=0.
\label{cqe}
\ee
Finally multiplying by $P^{ac}$ in the first of (\ref{cuarenta}) gives
\be
P^r_{\ s}(-P^a_{\ b}\bnb_aP^{sd}+p\bnb_bP^{sd})-P^{dr}E_b=0.
\label{zzz}
\ee
Multiplying this last equation by $P_z^{\ b}$ we obtain
$$
P^a_{\ z}P^r_{\ s}\bnb_aP^{sd}=0,
$$
which combined with (\ref{zzz}) implies 
\be
\bnb_bP^{rd}=\fr{1}{p}E_bP^{rd}.
\label{pre}
\ee
This last condition can be written with the indexes of the projector lowered if we use (\ref{cqe})
\be
\bnb_bP_{cd}=-\fr{1}{p}E_bP_{cd}.
\label{lll}
\ee
The dual conditions for the projectors $\Pi^{ab}$, $\Pi_{ab}$ coming from the vanishing of (\ref{def-N}) are
\be
\bnb_c\Pi^{ab}=\fr{1}{n-p}W_c\Pi^{ab},\ \bnb_c\Pi^a_{\ b}=0, \bnb_c\Pi_{ab}=-\fr{1}{p}W_c\Pi_{ab}.
\label{qrs}
\ee
Conversely, if equations (\ref{cqe}), (\ref{pre}), (\ref{lll}) and (\ref{qrs}) are assumed a simple calculation 
tells us that the tensors defined by (\ref{def-M}) and (\ref{def-N}) vanish. In fact all the above geometric 
conditions can be combined in a single simpler expression.
\begin{prop}
The following assertion is true
\bea
M_{abc}=\fr{1}{p}E_aP_{bc}-\fr{1}{n-p}W_a\Pi_{bc} 
\Longleftrightarrow\label{assertion0}\\
\bnb_aP_{bc}=-\fr{1}{p}E_aP_{bc},\ \bnb_a\Pi_{bc}=-\fr{1}{n-p}W_a\Pi_{bc}.
\label{assertion}
\eea
\label{ftp}
\end{prop} 

\begin{pf}
To prove this result we need the following
 formulae relating $\bnb_cP_{ab}$, $\bnb_cP^a_{\ b}$ and
$\nb_cP_{ab}$, $\nb_cP^a_{\ b}$. 
\bea
\bnb_aP_{bc}&=&\nb_aP_{bc}-\fr{1}{p}E_aP_{bc}-\fr{1}{2p}(E_bP_{ac}+E_cP_{ab})-\nonumber\\ 
&-&\fr{1}{2}(P_{cp}M^p_{\ ab}+P_{bp}M^p_{\ ac}),\label{identity-1}\\
2\bnb_aP^b_{\ c}&=&2\nb_aP^b_{\ c}+P^{bq}P^r_{\ c}M_{qra}-\Pi^{bq}P^r_{\ c}M_{qra}-
P^b_{\ q}M^q_{\ ac}+\fr{1}{n-p}W_c\Pi^b_{\ a}-\nonumber\\
&-&\fr{1}{p}E_cP^b_{\ a},
\label{identity-2}\\
\bnb_aP^{bc}&=&\nb_aP^{bc}+\fr{1}{p}E_aP^{bc}+\fr{1}{2(n-p)}(W^c\Pi^b_{\ a}+W^b\Pi^c_{\ a})-\nonumber\\
&-&\fr{1}{2}(M^b_{\ ar}P^{rc}+M^c_{\ ar}P^{rb}),
\label{identity-3}
\eea
There is a dual set of identities obtained through the replacements $P_{ab}\leftrightarrow\Pi_{ab}$, 
$E_a\leftrightarrow W_a$ and $p\leftrightarrow n-p$ (see paper I for a proof).
Now, if
conditions (\ref{assertion}) are true then expanding 
$\bnb_aP_{bc}$ by means of (\ref{identity-1}) we get
\bea
\nb_bP_{ac}&=&\fr{1}{2p}(E_aP_{bc}+E_cP_{ab})+\fr{1}{2}(P_{cp}M^p_{\ ba}+P_{ap}M^p_{\ bc})
\label{cdb-1}\\
\nb_b\Pi_{ac}&=&\fr{1}{2(n-p)}(W_a\Pi_{bc}+W_c\Pi_{ab})-\fr{1}{2}(\Pi_{cp}M^p_{\ ba}+\Pi_{ap}M^p_{\ bc}),
\label{cdb-2}
\eea
Substituting these expressions of the covariant derivatives of $P_{ab}$ and  
$\Pi_{ab}$ in the definition of $M_{abc}$ (equation (\ref{m-definition}))
yields
\be
P_{cp}M^p_{\ ab}=-\fr{1}{n-p}W_c\Pi_{ab},\ \ \Pi_{cp}M^p_{\ ab}=\fr{1}{p}E_cP_{ab},
\ee
whose addition leads to (\ref{assertion0}). 
Conversely, suppose that (\ref{assertion0}) holds. 
We may write this equation in the equivalent form
\be
\nb_bP_{ac}=\fr{1}{2p}(E_aP_{bc}+E_cP_{ab})-\fr{1}{2(n-p)}(W_a\Pi_{bc}+W_c\Pi_{ba}),
\label{25}
\ee
where the relation $M_{abc}=\nb_bP_{ac}-\nb_cP_{ab}-\nb_aP_{bc}$ has been again used.
Then inserting 
(\ref{25}) and (\ref{assertion0})
 into (\ref{identity-1}) gives us the condition on $\bnb_aP_{bc}$ at once.
 The calculation for $\bnb_a\Pi_{bc}$ is similar using identity (\ref{identity-1})
written in terms of $\Pi_{ab}$.\qed
\end{pf}
\begin{prop}
If (\ref{assertion0}) is true 
$\bnb_cP^a_{\ b}=0$
\label{up-down}
\end{prop}

\begin{pf}
To prove this we use identity 
(\ref{identity-2}) and replace $M_{abc}$ by (\ref{assertion0}) yielding
\be
\bnb_aP^b_{\ c}=\nb_aP^b_{\ c}-\fr{1}{2p}(E^bP_{ac}+E_cP^b_{\ a})+\fr{1}{2(n-p)}(W_c\Pi^b_{\ a}+W^b\Pi_{ac}),
\ee
which vanishes by (\ref{25}).\qed
\end{pf}

\begin{rem}
\label{remark}
{\em 
Note that in view of the above result condition (\ref{assertion0}) entails 
\be
\bnb_cP^{ab}=\fr{1}{p}E_cP^{ab},\ \ \bnb_c\Pi^{ab}=\fr{1}{n-p}W_c\Pi^{ab}.
\label{ftp-up}
\ee
Therefore all the conditions coming from (\ref{def-M}) and (\ref{def-N}) are summarized by
(\ref{assertion0}). This last equation can be written in the equivalent form $T_{abc}=0$
where
$$
T_{abc}\equiv M_{abc}-\fr{1}{p}E_aP_{bc}+\fr{1}{n-p}W_a\Pi_{bc} 
$$
As explained in paper I 
the tensor $T_{abc}$ plays a key role in the 
geometric characterization of {\em conformally 
separable pseudo-Riemannian manifolds}. There we proved that a pseudo-Riemannian manifold is locally 
conformally separable with the tensors $P_{ab}$ and $\Pi_{ab}$ as the leaf metrics (see definition \ref{separable}) if 
and only if $T_{abc}=0$ which means that the complete integrability conditions obtained so far 
bear a clear geometrical meaning.}
\end{rem}

Some of the first integrability conditions achieve a great simplification if $T_{abc}$ vanishes. For instance 
(\ref{cond-ligaduras-down})-(\ref{mixed-invariance}) become zero identically under this condition as is obvious from 
propositions \ref{ftp} and \ref{up-down} so we do not need to care about these integrability
conditions any more. Equation
(\ref{cond-g1}) acquires the invariance law
\be
\lie T^d_{\ cab}=0.
\label{cond-g11}
\ee
Other simplifications will be shown in the forthcoming analysis.

\medskip
\noindent
{\em Equations (\ref{cond-d}) and (\ref{cond-e}).} If $T_{abc}=0$ then $\La^a_{\ bc}=0$, $\bLa^a_{\ bc}=0$ 
so these equations take the form 
\bea
-\bphi_qP^q_{\ r}T^r_{\ ceb}=2\lie\bnb_{[e}L^0_{b]c}+\bchi_rE^r_{\ ceb},\label{cond-df}\\ 
-\bchi_q\Pi^q_{\ r}T^r_{\ ceb}=2\lie\bnb_{[e}L^1_{b]c}+\bphi_rF^r_{\ ceb},\label{cond-ef}  
\eea
where
\bea
E^r_{\ ceb}=\fr{2}{2-p}(\Pi^q_{\ c}\Pi^r_{\ [e}L^0_{\ b]q}+\Pi^r_{\ c}\Pi^q_{\ [e}L^0_{b]q}-\Pi^{rq}\Pi_{c[e}L^0_{b]q})\\
F^r_{\ ceb}=\fr{2}{2-n+p}(P^q_{\ c}P^r_{\ [e}L^1_{\ b]q}+P^r_{\ c}P^q_{\ [e}L^1_{b]q}-P^{rq}P_{c[e}L^1_{b]q}).
\eea
If  (\ref{cond-df}) and (\ref{cond-ef}) are to be true for every value of 
$\bphi_q$ and $\bchi_q$ we find the conditions of complete integrability
$$
P^q_{\ r}T^r_{\ ceb}=0,\ \Pi^q_{\ r}T^r_{\ ceb}=0\Longrightarrow T^q_{\ ceb}=0,
$$
so equation (\ref{cond-g11}) is trivially fulfilled.
To proceed further with the calculations we need a lemma 
\begin{lem}
If $\bnb_aP^b_{\ c}=0$ then 
$$
L^0_{bq}\Pi^q_{\ c}=0,\ \ L^1_{bq}P^q_{\ c}=0.
$$
\end{lem}

\begin{pf}
These properties are proven through the Ricci identity applied to the tensors $P^b_{\ c}$, $\Pi^b_{\ c}$ which 
take a remarkably simple form under our conditions (we only perform the calculations for the tensor $P^a_{\ b}$)
\be
0=\bnb_e\bnb_aP^b_{\ c}-\bnb_a\bnb_eP^b_{\ c}=\bR^b_{\ qea}P^q_{\ c}-\bR^q_{\ cea}P^b_{\ q},
\label{ricci-p}
\ee 
whence
$$
L^0_{bq}\Pi^q_{\ c}=2P^d_{\ r}\bR^r_{\ qdb}\Pi^q_{\ c}+\fr{2}{p}(\Pi^d_{\ b}P^r_{\ q}\bR^q_{\ rds}\Pi^s_{\ c}).
$$
The first term of this expression is zero according to (\ref{ricci-p}) and the second one can be transformed 
by means of the first Bianchi identity into
$$
\Pi^d_{\ b}\Pi^s_{\ c}P^r_{\ q}\bR^q_{\ rds}=-\Pi^d_{\ b}\Pi^s_{\ c}(P^r_{\ q}\bR^q_{\ dsr}+P^r_{\ q}\bR^q_{\ srd})=
-\Pi^d_{\ b}\Pi^s_{\ c}(P^q_{\ d}\bR^r_{\ qsr}+P^q_{\ s}\bR^r_{\ qrd}),
$$
which also vanishes.\qed
\end{pf}
Therefore conditions (\ref{cond-df}) and (\ref{cond-ef}) are further simplified to 
\be
\lie\bnb_{[e}L^0_{b]c}=0,\ \ \ 
\lie\bnb_{[e}L^1_{b]c}=0.
\label{further}
\ee
It is our next aim to show that indeed these two equations are identities if $T^d_{\ cab}=0$.
\begin{lem}
If $p\neq 3$, $n-p\neq 3$ and $T_{abc}=0$ then
$$
T^d_{\ cab}=0\Longrightarrow \bnb_{[e}L^0_{b]c}=0,\ \ \bnb_{[e}L^1_{b]c}=0
$$
\label{long-lema}
\end{lem}
\begin{pf}
To prove this we start from the identity
\bea
& &(2-p)[\bnb_e(P^e_{\ q}T^q_{\ cab})+\bnb_b(P^d_{\ q}T^q_{\ cda})-\bnb_a(P^d_{\ q}T^q_{\ cdb})]=\nonumber\\
&=&2P^e_{\ c}\bnb_{e}L^0_{[ba]}
-2P^d_{\ c}\bnb_{[b|}L^0_{d|a]}+2p\bnb_{[b}L^0_{a]c}+
2P^d_{\ [b}\bnb_{a]}L^0_{dc}-2P^e_{\ [b|}\bnb_eL^0_{|a]c}+\nonumber\\
&+&2P^{qe}\bnb_eL^0_{[a|q}P_{|b]c}+2P^{qd}P_{c[a}\bnb_{b]}L^0_{dq},
\label{ident-1} 
\eea
which is  obtained from equation (\ref{tensor-tt}) and the second Bianchi identity for the tensor 
$\bR^a_{\ bcd}$ if the condition $T_{abc}=0$ holds. 
By assumption 
the left hand side of this identity vanishes so we only need to study the 
right hand side equated to zero. Multiplying such equation by $P^{ca}$ we get
\be
\hspace{-1cm}
2(1-p)P^{ae}\bnb_eL^0_{ba}-2(1-p)P^{da}\bnb_bL^0_{da}+2P^d_{\ b}P^{ac}\bnb_aL^0_{dc}-2P^e_{\ b}P^{ac}\bnb_eL^0_{ac}=0,
\label{two-terms}
\ee
and a further multiplication by $P_z^{\ b}$ yields
$$
P_z^{\ r}P^{qe}\bnb_eL^0_{rq}=P_z^{\ r}P^{qd}\bnb_rP_{qd}.
$$
This last property implies that the last two terms of (\ref{two-terms}) are zero
from which we conclude that
$$
P^{qe}\bnb_eL^0_{rq}=P^{qd}\bnb_rP_{qd},
$$
and hence (\ref{ident-1}) becomes
\be
2P^e_{\ c}\bnb_{e}L^0_{[ba]}
-2P^d_{\ c}\bnb_{[b|}L^0_{d|a]}
+2p\bnb_{[b}L^0_{a]c}+2P^d_{\ [b}\bnb_{a]}L^0_{dc}=0.
\label{ident-2}
\ee
 If we multiply this last equation by $P_z^{\ a}P_t^{\ b}$  we obtain
\be
S_{ctz}-S_{tcz}+S_{zct}-S_{czt}+(p-2)(S_{tzc}-S_{ztc})=0,
\label{Sabc}
\ee
where
$$
S_{ztc}\equiv P_z^{\ r}P_t^{\ s}\bnb_rL^0_{sc}.
$$
Permuting indexes in (\ref{Sabc}) we get the equations
\be
\left.
\begin{array}{c}
S_{ctz}-S_{tcz}+S_{zct}-S_{czt}+(p-2)(S_{tzc}-S_{ztc})=0\\
S_{ztc}-S_{tzc}+S_{czt}-S_{zct}+(p-2)(S_{tcz}-S_{ctz})=0\\
S_{ctz}-S_{tcz}+S_{tzc}-S_{ztc}+(p-2)(S_{zct}-S_{czt})=0
\end{array}\right\}.
\ee
Setting the variables $x=S_{ctz}-S_{tcz}$, $y=S_{zct}-S_{czt}$, $w=S_{tzc}-S_{ztc}$ we deduce
that previous equations form a homogeneous system in these variables whose matrix is
$$
\left(\begin{array}{ccc}
1&1&p-2\\
2-p&-1&-1\\
1&p-2&1\end{array}\right)\Rightarrow \left|
\begin{array}{ccc}
1&1&p-2\\
2-p&-1&-1\\
1&p-2&1\end{array}\right|=-p(p-3)^2.
$$
So unless $p=3$ ($p=0$ makes no sense in the current context) we conclude that $x=y=w=0$ and hence 
$$
P_{[a}^{\ r}P_{b]}^{\ s}\bnb_rL^0_{sc}=0.
$$
Application of this in the expression resulting of multiplying (\ref{ident-2}) by $P^r_{\ b}$ leads to
$$
P^d_{\ c}\bnb_aL^0_{db}-P^e_{\ c}\bnb_eL^0_{ab}+(p-1)(P^r_{\ b}\bnb_rL^0_{ac}-P^r_{\ b}\bnb_aL^0_{rc})=0.
$$
By setting $Q_{acb}\equiv P^d_{\ c}\bnb_aL^0_{db}-P^e_{\ c}\bnb_eL^0_{ab}$ we can rewrite this as
$$
Q_{acb}-(p-1)Q_{abc}=0,\Rightarrow Q_{abc}-(p-1)Q_{acb}=0,
$$
which entails $Q_{abc}=0$ (recall that $p\neq 1$ by definition of $L^0_{ab}$). 
This last property applied to (\ref{ident-2}) yields 
$$
\bnb_{[b}L^0_{a]c}=0,
$$ 
as desired. The result for $L^1_{ab}$ is proven in a similar way.\qed
\end{pf}

\medskip
\noindent
{\em Equations (\ref{cond-a1}) and (\ref{cond-a2})}. The analysis of these conditions is performed 
by means of the following result.
\begin{prop}
If $T_{abc}=0$ then 
\bnr
\fr{1}{2-p}(L^0_{ab}-L^0_{ba})-\fr{1}{p}(\bnb_aE_b-\bnb_bE_a)=0,\\ 
\fr{1}{2-n+p}(L^1_{ab}-L^1_{ba})-\fr{1}{n-p}(\bnb_aW_b-\bnb_bW_a)=0.
\enr
\end{prop}
\begin{pf}
We only carry on the proof for the first identity as the calculations are similar for the second one. We start from the identity
$$
E_a=-P^{bc}\bnb_aP_{bc},
$$ 
which is easily obtained from (\ref{identity-1}). Using 
this we may write
\be
\hspace{-1cm}\bnb_aE_b-\bnb_bE_a=-P^{qr}(\bnb_a\bnb_bP_{qr}-\bnb_b\bnb_aP_{qr})-\bnb_aP^{qr}\bnb_bP_{qr}+\bnb_bP^{qr}\bnb_aP_{qr}.
\label{qxi}
\ee
The expression in brackets can be transformed by the Ricci identity into 
$2P^q_{\ r}\bR^r_{\ qab}$.
If we impose now the condition $T_{abc}=0$, then combination of proposition 
\ref{ftp} and remark \ref{remark} entails
$$
\bnb_aP^{qr}\bnb_bP_{qr}=-\fr{1}{p}E_aE_b=\bnb_bP^{qr}\bnb_aP_{qr},
$$
Therefore after these manipulations
equation (\ref{qxi}) yields
\be
\bnb_aE_b-\bnb_bE_a=2P^q_{\ r}\bR^r_{\ qab}.
\label{dos}
\ee
On the other hand from 
(\ref{L0}) and applying the first Bianchi identity it is easy to obtain
\be
L^0_{ab}-L^0_{ba}=\fr{2(2-p)}{p}P^r_{\ q}\bR^q_{\ rab}.
\label{uno}
\ee
Combination of (\ref{uno}) and (\ref{dos}) leads to the desired result.\qed
\end{pf}

\medskip
\noindent
{\em Equations (\ref{cond-b}) and (\ref{cond-c})}
\begin{prop}
If $T^d_{\ cab}=0$ and $T_{abc}=0$ then (\ref{cond-b}) and (\ref{cond-c}) are 
identities.
\end{prop}
\begin{pf}
The left hand side of both equations vanishes trivially if $T^d_{\ cab}=0$ so we just need to show that 
the right hand side vanishes as well. The characterization of the condition $T_{abc}=0$ in terms of the 
covariant derivatives of the projectors 
(equations (\ref{assertion}), (\ref{ftp-up})) entails
$$
\Up^{sc}_{\ b}=\fr{2-p}{p}E_bP^{sc},\ \ \bUp^{sc}_{\ b}=\fr{2-n+p}{n-p}W_b\Pi^{sc}, 
$$
which means that the terms of (\ref{cond-b}) and (\ref{cond-c}) containing these tensors are zero.
The property $\bnb_cP^a_{\ b}=\bnb_c\Pi^a_{\ b}=0$ can be used now to get rid of some terms and simplify others
on these couple of equations getting
\bnr
0=\bchi_a\bnb_{[c}E_{b]}+\bchi_{[b}\bnb_{c]}E_a+\bchi_a\bnb_{[b}E_{c]}+\bchi_{[c}\bnb_{b]}E_a\\
0=\bphi_a\bnb_{[c}W_{b]}+\bphi_{[b}\bnb_{c]}W_a+\bphi_a\bnb_{[b}W_{c]}+\bphi_{[c}\bnb_{b]}W_a,
\enr
which are obviously identities .\qed
\end{pf}

\subsection*{Constraints}
If $T_{abc}=0$ then 
using (\ref{assertion}) and proposition \ref{up-down}
equation (\ref{cond-ligaduras-down}) becomes 
$$
-\fr{1}{p}\lie(E_cP_{ab})=-\fr{\phi}{p} E_cP_{ab}+\sphi_cP_{ab},\ 
-\fr{1}{n-p}\lie(W_c\Pi_{ab})=\fr{-\chi}{n-p}W_c\Pi_{ab}+\schi_c\Pi_{ab},
$$
which is easily seen to be an identity if we apply 
(\ref{ligaduras})-$B$.

All our calculations are thus summarized in the next result 
which is one of the most important of this paper.
\begin{thm}[Complete integrability conditions]
The first integrability conditions calculated for bi-conformal 
vector fields are identically fulfilled in the neighbourhood of 
a point if and only if in such neighbourhood  
\begin{center}
\fbox{$T_{abc}=0,\ \ \ T^d_{\ cab}=0$},
\end{center}
\label{complete-integrability}
whenever $P^a_{\ a}$, $\Pi^a_{\ a}\neq 3$. 
\label{integrability-complete}
\end{thm}

\section{Geometric characterization of bi-conformally flat 
pseudo-Rieman\-ni\-an manifolds}
\label{bi-flat}

Once we have found the mathematical characterization of the spaces admitting a maximum number of bi-conformal vector 
fields we must next settle if there is actually any space whose metric tensor complies with the conditions stated in theorem
\ref{complete-integrability} or on the contrary there are no pseudo-Riemannian manifolds fulfilling such requirement.
Indeed, we will find that each geometric condition has a separate meaning related to the geometric characterization of 
certain {\em separable} pseudo-Riemannian manifolds. Hence the tensors $T_{abc}$ and $T^a_{\ bcd}$ bear a geometric 
interest on their own regardless of the existence of bi-conformal vector fields on the pseudo-Riemannian manifold where they are 
defined. Before addressing this issue 
we need some preliminary definitions.
\begin{defn}
The pseudo-Riemannian manifold $(V,\rmg_{ab})$ is said to be conformally separable at the point $q\in V$ if there exists a 
local coordinate chart 
$x\equiv\{x^1,\dots,x^n\}$ based at $q$ in which the metric tensor takes the form
\be
\rmg_{ab}(x)=\left\{
\begin{array}{c}
\Xi_1(x)G_{\alpha\beta}(x^\g),\ 1\leq\alpha,\beta,\g\leq p\\
\Xi_2(x)G_{AB}(x^C),\ p+1\leq A,B,C\leq n\\
0\ \ \mbox{otherwise.}
\end{array}\right.
\label{metric-separable}
\ee
where $\Xi_1$, $\Xi_2$ are $C^2$ functions on the open set defining the coordinate chart.
$(V,\rmg_{ab})$ is conformally separable if it is so at every point $q\in V$. Any of the tensors 
$\Xi_1G_{\alpha\beta}$, $\Xi_2G_{AB}$ shall be called leaf metric.
\label{separable}
\end{defn}
In paper I we proved that a pseudo-Riemannian manifold is conformally separable at a point $q$ if 
and only if the tensor $T_{abc}$ is zero in a neighbourhood of $q$. In this case the orthogonal 
projectors $P_{ab}$ and $\Pi_{ab}$ used to construct $T_{abc}$
generate integrable distributions $D$ and $D'$ (see considerations coming 
after definition \ref{Bi-conformal}) and hence the set of integral manifolds 
of each distribution is a foliation in a neighbourhood of $q$. These foliations 
allow us to construct a local coordinate system around $p$ in which the metric tensor
takes the form (\ref{metric-separable}) and the only non-vanishing
components of the projectors are
$$
P_{\alpha\beta}=\Xi_1(x)G_{\alpha\beta}(x^{\g}),\ 
\Pi_{AB}=\Xi_2(x)G_{AB}(x^{C}),
$$
from which we conclude that $P_{ab}$ and $\Pi_{ab}$ are just 
the leaf metrics. 
The advantage of the characterization of conformally separable pseudo-Riemannian
manifolds given in paper I is clear as it is invariant and coordinate-free. Therefore
we conclude that any pseudo-Riemannian manifold fulfilling the conditions of theorem 
\ref{integrability-complete} must be conformally separable.

In what is to follow we will work with conformally separable pseudo-Riemannian 
manifolds at a point and thus we will identify  
$P_{ab}$, $\Pi_{ab}$ and the leaf metrics.

Bi-conformally flat spaces are a particular and interesting case of conformally 
se\-pa\-rable pseudo-Riemannian 
manifolds defined by the requirement that the metrics $G_{\alpha\beta}$, 
$G_{AB}$ be flat (so the leaf metrics
are conformally flat). To our knowledge the study of bi-conformally flat spaces 
has never been tackled 
in the literature as opposed to many other conformally separable pseudo-Riemannian 
manifolds.
Next, we fill up this gap proving a local characterization of bi-conformally flat pseudo-Riemannian 
manifolds along the same lines as in the conformally separable case.  
\begin{thm}
A conformally separable pseudo-Riemannian manifold at a point 
$q$ with leaf metrics of rank greater than 3 
is bi-con\-for\-ma\-lly flat at the same point if and only if the tensor
 $T^a_{\ bcd}$ constructed from its leaf metrics is identically zero in 
a neighbourhood of $q$.
\label{bi-conformal-char} 
\end{thm}

\begin{pf} 
We choose the same coordinates and notation for our conformally separable space
 as in definition \ref{separable} (from now on we will label tensor 
indexes with lowercase Greek and uppercase Latin letters according to the leaf of the 
separation they refer to). In the local coordinates of (\ref{metric-separable}) the non-vanishing Christoffel symbols 
are 
\bnr
\Chr^{\alpha}_{\ \beta\g}=\fr{1}{2\Xi_1}G^{\alpha\rho}(\d_{\beta}(\Xi_1G_{\alpha\rho})+
\d_{\g}(\Xi_1G_{\rho\beta})-\d_{\rho}(\Xi_1G_{\beta\g})),\\
\Chr^{A}_{\ BC}=\fr{1}{2\Xi_1}G^{AD}(\d_{B}(\Xi_1G_{CD})+
\d_{C}(\Xi_1G_{DB})-\d_{D}(\Xi_1G_{BC})),\\
\Chr^{\alpha}_{\beta A}=\fr{1}{2\Xi_1}\delta^{\alpha}_{\ \beta}\d_A\Xi_1,\ \ 
\Chr^A_{\ B\alpha}=\fr{1}{2\Xi_2}\delta^A_{\ B}\d_{\alpha}\Xi_2,
\enr
where, 
$$
G^{\alpha\beta}G_{\beta\rho}=\delta^{\alpha}_{\ \rho},\ 
G^{AC}G_{CB}=\delta^A_{\ B},\ 
$$
and the tensors $G^{\alpha\beta}$ and $G^{AB}$ are used to raise 
Greek and Latin indexes respectively. 
The components of $M_{abc}$, $E_a$, $W_a$ are 
(henceforth, the components not shown in 
an explicit tensor representation are understood to be zero) 
\bea
M_{\alpha AB}=\d_{\alpha}(\Xi_2G_{AB}),\ \ 
M_{A\alpha\beta}=-\d_{A}(\Xi_1G_{\alpha\beta}),\nonumber\\
E_{A}=-\d_{A}\log|\mbox{det}(\Xi_1G_{\alpha\beta})|,\ \ 
W_{\alpha}=-\d_{\alpha}\log|\mbox{det}(\Xi_2G_{AB})|.
\label{27}
\eea
Therefore we get for the components of the bi-conformal connection
\bea
\bar{\Chr}^{\alpha}_{\ \beta\phi}&=&\fr{1}{2\Xi_1}(\delta^{\alpha}_{\ \beta}\d_{\phi}\Xi_1+\delta^{\alpha}_{\ \phi}\d_{\beta}\Xi_1-
G^{\alpha\rho}G_{\beta\phi}\d_{\rho}\Xi_1)+\Chr^{\alpha}_{\ \beta\phi}(G),\nonumber\\
\bar{\Chr}^A_{\ BC}&=&\fr{1}{2\Xi_2}(\delta^{A}_{\ B}\d_{C}\Xi_2+\delta^{A}_{\ C}\d_{B}\Xi_2-
G^{AR}G_{BC}\d_{R}\Xi_2)+\Chr^{A}_{\ BC}(G)\nonumber\\
\bar{\Chr}^{\alpha}_{\ \beta C}&=&\bar{\Chr}^{A}_{\ B\phi}=0,
\label{new-components}
\eea
where $\Chr^{\alpha}_{\ \beta\phi}(G)$ and $\Chr^{A}_{\ BC}(G)$ are the Christoffel symbols of the metrics 
$G_{\alpha\beta}$ and $G_{AB}$ respectively.
Using this we deduce after some algebra
\bea
\bR^{\alpha}_{\ \beta\phi\g}&=&R^{\alpha}_{\ \beta\phi\g},\ \bR^{A}_{\ BCD}=R^A_{\ BCD},\ 
\bR^{\alpha}_{\ \beta F\phi}=\d_F\bar{\Chr}^{\alpha}_{\ \phi\beta},\ \bR^A_{\ BF\phi}=-\d_{\phi}\bar{\Chr}^A_{\ FB}
\nonumber\\
L^0_{\alpha\beta}&=&2R_{\alpha\beta}+\fr{R^{\g}_{\ \g}}{1-p}\rmg_{\alpha\beta},\ \ 
L^1_{AB}=2R_{AB}+\fr{R^C_{\ C}}{1-n+p}\rmg_{AB},\ L^0_{\alpha A}=L^1_{A\alpha}=0,\nonumber\\
L^0_{\ A\alpha}&=&\fr{2(2-p)}{p}\d_A\bar{\Chr}^{\rho}_{\ \alpha\rho},\ L^1_{\alpha A}=\fr{2(2-n+p)}{n-p}\d_{\alpha}\bar{\Chr}^{R}_{\ AR}
\label{cua} 
\eea
where $R_{\alpha\beta}$, $R_{AB}$, $R^{\g}_{\ \g}$ and $R^{C}_{\ C}$ are the Ricci tensors and Ricci scalars of 
$R^{\alpha}_{\ \beta\phi\g}$ and $R^A_{\ BCD}$ respectively (these 
curvature tensors are calculated from the leaf metrics and not 
from their conformal counterparts $G_{\alpha\beta}$ and $G_{AB}$). 
From this and equations (\ref{tensor-tt}), (\ref{new-components}) we get 
that the only non-vanishing components of the tensor $T^a_{\ bcd}$ are
\be
T^{\alpha}_{\ \beta\phi\g}=2C^{\alpha}_{\ \beta\phi\g},\ \ \ T^A_{\ BCD}=2C^A_{\ BCD},
\label{weyl}
\ee
being $C^{\alpha}_{\ \beta\phi\g}$ and $C^A_{\ BCD}$ the Weyl tensors constructed from each leaf metric 
through the relations
\bnr
C^{\alpha}_{\ \beta\phi\g}=R^{\alpha}_{\ \beta\phi\g}+\fr{1}{2-p}(\rmg_{\beta[\g}L^0_{\phi]\rho}\rmg^{\rho\alpha}
+\delta^{\alpha}_{\ [\phi}L^0_{\g]\beta}),\\ 
C^A_{\ BCD}=R^A_{\ BCD}+\fr{1}{2-n+p}(\rmg_{B[D}L^1_{C]E}\rmg^{EA}+\delta^A_{\ [C}L^1_{D]B}).
\enr
Hence $T^{\alpha}_{\ \beta\phi\g}$ and  $T^A_{\ BCD}$ are both zero if and only if the leaf metrics are both conformally flat
which proves the theorem.\qed 
\end{pf}

\begin{rem}
{\em Note that the tensors $\bR^{\alpha}_{\ \alpha F\phi}$ and $\bR^{A}_{\ AF\phi}$ (no 
summation) do not vanish in general 
which means that the bi-conformal connection does not stem from a metric tensor in this case.}
\end{rem}
From the calculations performed above it is clear that theorem \ref{bi-conformal-char} 
can be generalized to 
conformally separable spaces in which only one of the leaf metrics 
is conformally flat. To that end, we define the tensor 
\be
T(P)^a_{\ bcd}\equiv P^a_{\ r}P^q_{\ b}P^{s}_{\ c}P^t_{\ d}T^r_{\ qst}.
\label{t-p}
\ee
\begin{thm}
Under the assumptions of theorem \ref{bi-conformal-char} a leaf metric is conformally flat if and only if 
the tensor $T(P)^a_{\ bcd}$ calculated from the leaf metric $P_{ab}$ 
is equal to zero.  
\label{bi-conformal-char-2}
\end{thm}
\begin{pf}
From the proof of theorem \ref{bi-conformal-char}  and (\ref{t-p}) 
we deduce that for a conformally separable 
pseudo-Riemannian manifold the only non-vanishing components of $T(P)^a_{\ bcd}$ are 
$$
T(P)^{\alpha}_{\ \beta\g\phi}=T^{\alpha}_{\ \beta\g\phi}=2C^{\alpha}_{\ \beta\g\phi},\ \ 
$$
so the vanishing of $T(P)^a_{\ bcd}$ implies that the Weyl tensor calculated from the corresponding leaf metric is zero 
as well.\qed
\end{pf}

In the case of any of the leaf metrics being of rank three 
it is clear from the above 
that the corresponding tensor 
$T(P)^a_{\ bcd}$ will be zero as the Weyl tensor of any three dimensional 
pseudo-Riemannian metric 
vanishes identically. Hence the results presented so far cannot be used to characterize 
conformally separable pseudo-Riemannian manifolds with conformally flat leaf metrics.
 This lacking is remedied in the next theorem.
\begin{thm}
A conformally separable pseudo-Riemannian manifold at a point $q$ 
has a conformally flat leaf metric of rank three if and only
if the condition
\be
\bnb_{[a}L^0_{b]c}=0, 
\label{rank-3}
\ee
holds in a neighbourhood of $q$
where $L^0_{ab}$ is calculated from the leaf metric
 $P_{ab}$ by means of (\ref{L0}).  
\label{case-3}
\end{thm}   
\begin{pf} 
To show this we will rely on the notation and calculations performed in the proof of theorem 
\ref{bi-conformal-char}. 
If the manifold is conformally separable,
then the components of the tensor $\bnb_aL^0_{bc}$ are
\bnr
\bnb_{\alpha}L^0_{\beta\epsilon}&=&\nb_{\alpha}L^0_{\beta\epsilon},\ \ 
\bnb_{\alpha}L^0_{B\epsilon}=\fr{2(2-p)}{p}(\d^2_{\alpha B}\bar{\Chr}^{\rho}_{\ \epsilon\rho}-
\bar{\Chr}^{\delta}_{\ \alpha\epsilon}\d_B\bar{\Chr}^{\rho}_{\delta\rho})
,\ \bnb_AL^0_{\beta\epsilon}=\d_AL^0_{\beta\epsilon},\\
\bnb_AL^0_{B\epsilon}&=&\d_AL^0_{B\epsilon}-\bar{\Chr}^C_{AB}L^0_{C\epsilon},
\enr
where $\nb_{\alpha}$ is the covariant derivative compatible with the 
metric $\Xi_1(x)G_{\alpha\beta}$. Trivially
\be
\bnb_{[\alpha}L^0_{\beta]\epsilon}=\nb_{[\alpha}L^0_{\beta]\epsilon},\ \
\label{eq}
\ee
here the tensor $\nb_{[\alpha}L^0_{\beta]\epsilon}$ is the Schouten tensor of the 3-metric $\Xi_1(x)G_{\alpha\beta}$ 
which vanish if and only if $G_{\alpha\beta}$ is flat. 
Therefore to finish the proof of this theorem 
we must show that all the remaining components of 
$\bnb_{[a}L^0_{b]c}$ are zero. These are
$$
\bnb_{[A}L^0_{B]\epsilon}=\fr{2(2-p)}{p}\d^2_{[AB]}\bar{\Chr}^{\rho}_{\ \epsilon\rho},\ 
2\bnb_{[\alpha}L^0_{B]\epsilon}=\fr{2(2-p)}{p}(\d^2_{\alpha B}\bar{\Chr}^{\rho}_{\ \epsilon\rho}-
\bar{\Chr}^{\delta}_{\ \alpha\epsilon}\d_B\bar{\Chr}^{\rho}_{\delta\rho})-\d_BL^0_{\alpha\epsilon}.
$$
Clearly the first expression is zero and the second one is worked out by replacing 
the connection coefficients by their expressions given in (\ref{new-components}) and 
using the identity 
\bnr
L^0_{\alpha\beta}&=&L^0_{\alpha\beta}(G)+(2-p)(2\sigma_{\alpha\beta}+G_{\alpha\beta}(\d\sigma)^2),\\ 
\sigma&=&\fr{1}{2}\log|\Xi_1(x)|,\ \sigma_{\alpha\beta}=\hat{\nb}_{\alpha}\hat{\nb}_{\beta}\sigma
-\d_{\alpha}\sigma\d_{\beta}\sigma,\ 
(\d\sigma)^2=G^{\alpha\beta}\d_{\alpha}\sigma\d_{\beta}\sigma,
\enr
where $L^0_{\alpha\beta}(G)$ is calculated using the curvature tensors computed from $G_{\alpha\beta}$ and 
$\hat{\nb}$ is the connection compatible with this metric. The sought result comes after some 
simple algebraic manipulations.
\qed
\end{pf} 

All in all the geometric conditions proven in theorems 
\ref{bi-conformal-char}, \ref{bi-conformal-char-2} and \ref{case-3} provide 
a set of equations which can be used to determine if a given 
pseudo-Riemannian manifold $(V,\rmg_{ab})$ is  
conformally separable at a point with respect to 
a pair of leaf metrics $P_{ab}$, $\Pi_{ab}$ and to decide 
if any of the leaf metrics is conformally flat. 
Note that all our characterizations are coordinate-free and 
in fact they can be used as intrinsic definitions of conformal
separability and bi-conformal flatness at a point as opposed to 
definition (\ref{separable}) which is coordinate-dependent. What is more 
one can study specific examples of pseudo-Riemannian manifolds in 
which tensors $P_{ab}$ and $\Pi_{ab}$ with the properties (\ref{properties})
are defined
and check if the afore-mentioned conditions are fulfilled which is far more 
easier than trying to prove the existence of the coordinate system of 
(\ref{metric-separable}) (see example \ref{example1} for a practical
application). 

\subsection{Bi-conformally flat spaces as spaces with a maximal 
number of bi-conformal vector fields}
\label{bi-maximal}
From the above it is clear that the conditions imposed 
by theorem \ref{complete-integrability} are satisfied by a nontrivial 
set of pseudo-Riemannian manifolds. Thus we deduce that 
a bi-conformally flat pseudo-Riemannian manifold 
can be also characterized by the existence of a maximum number 
of bi-conformal vector fields. On the other hand it is 
straightforward to check (proposition 6.1 of \cite{BI-CONFORMAL}) 
that for these spaces any conformal Killing vector of the leaf metrics 
is a bi-conformal vector field of the metric $\rmg_{ab}$. As the number 
of conformal Killing vectors for each leaf metric is the biggest possible 
as well we get at once that for a bi-conformally flat space 
the total number of linearly independent bi-conformal vector fields is 
$$
\fr{1}{2}(p+1)(p+2)+\fr{1}{2}(n-p+1)(n-p+2),\ \ p,\ n-p\neq 2
$$       
which is the upper bound $N$ for the dimension of any finite dimensional Lie 
algebra of bi-conformal 
vector fields. Summing up we obtain the following result
\begin{thm}
A pseudo-Riemannian manifold possesses $N$ linearly independent 
bi-conformal vector fields ($P^a_{\ a}> 3$, $\Pi^a_{\ a}> n-3$) 
if and only if it is bi-conformally flat.\qed
\label{N}
\end{thm}
Bi-conformally flat spaces in which any of the leaf metrics has rank three 
also admit $N$ linearly independent bi-conformal 
vector fields (in fact the complete integrability conditions are satisfied 
for these spaces as a result 
of theorem \ref{case-3}). However, we do not know yet if there are spaces with 
$N$ linearly independent bi-conformal vector fields with either of the projectors 
$P_{ab}$ or $\Pi_{ab}$ projecting on a 3-dimensional vector 
space other than bi-conformally flat spaces. This is so because 
in such case the complete integrability conditions 
(\ref{further}) may in principle be fulfilled by other conformally 
separable spaces not necessarily with conformally 
flat leaf metrics. The true extent of this assertion and the 
complete characterization of spaces with $N$ linearly 
independent bi-conformal vector fields under these circumstances 
will be placed elsewhere.

\section{Examples}
We present here examples illustrating how our techniques 
work in practical cases. All the algebraic calculations can be 
performed with any of the computer algebra systems available today 
(the system used here was GRTensorII \cite{GRTENSOR}).

\label{examples}
\begin{exmp}
\label{example1}
\em Consider the four dimensional pseudo-Riemannian 
metric given by 
$$
ds^2=\Psi^2\sin^2\z(dt+d\phi)^2-\alpha dt^2+B^2(dr^2+r^2d\z^2)
$$
where $-\infty<t<\infty$, $0<r<\infty$, $0<\z<\pi$, $0<\phi<2\pi$ and  
$$
\Psi=\Psi(r,\z),\ \Phi=\Phi(r,\z),\ 
B=B(r,\z),\ \alpha=\alpha(r,\z).
$$
Let us prove that this metric is 
conformally separable with respect to the foliation 
defined by the condition $t=const$. To that end we need to find out 
the orthogonal projector  $P_{ab}$ whose associated distribution $D$
generates the above foliation. The distribution $D$ is clearly defined 
by the set of vector fields $\{\d/\d r,\d/\d\theta,\d/\d\phi\}$ 
which are already orthogonal and hence equation (\ref{distribution}) 
tells us how to construct $P_{ab}$ and $\Pi_{ab}$. The only non-vanishing 
components of $P_{ab}$ are
$$
P_{tt}=\Psi^2\sin^2\z,\ P_{rr}=B^2,\ P_{\z\z}=r^2B^2,\ 
P_{\phi\phi}=\Psi^2\sin^2\z,\ 
P_{t\phi}=\Psi^2\sin^2\z.
$$
and from this we can easily check that the tensor $T_{abc}$ is zero 
identically. The tensor $\bnb_{[a}L^0_{b]c}$ is not zero unless further 
restrictions are imposed. For instance if we set $\Psi^2=r^2(1+r^2)$, 
$B^2=1+r^2$, $\alpha=-r^2$ the metric tensor becomes
\bea
ds^2=r^2(1+(1+r^2)\sin^2\z)dt^2+2r^2(1+r^2)\sin^2\z dtd\phi+\nonumber\\
+(1+r^2)(dr^2+r^2d\z^2+r^2\sin^2\z d\phi^2),
\label{case}
\eea
 and a calculation shows that condition (\ref{rank-3}) holds. 
Hence according to theorem \ref{case-3} the metric is locally 
bi-conformally flat. Note that this is by no means evident in 
the coordinate system of (\ref{case}) and so our method has a clear 
advantage. In this case we can go even further and find the explicit 
coordinates bringing (\ref{case}) into the canonical form 
(\ref{metric-separable}) which is  
$$
ds^2=(x^2+y^2+z^2)dT^2+(1+x^2+y^2+z^2)(dx^2+dy^2+dz^2),
$$
being the coordinate change 
$$
T=t,\ x=r\sin\z\cos(t+\phi),\ y=r\sin\z\sin(t+\phi),\ z=r\cos\z.
$$     
\end{exmp}

\begin{exmp}
\label{example2}
\em
In the foregoing results we have only concentrated on conformally separable 
pseudo-Riemannian manifolds but nothing was said about manifolds 
foliated by conformally flat hypersurfaces and not 
conformally separable. 
To illustrate this case let us consider the four dimensional 
pseudo-Riemannian manifold $V$ 
covered by a single coordinate chart $x=\{x^1,x^2,x^3,x^4\}$ and whose 
metric tensor is
\be
\hspace{-.4cm}ds^2=\Phi(x)[(dx^1)^2+(dx^2)^2+(dx^3)^2]+2\sum_{i=1}^3\beta_i(x)dx^idx^4+
\Psi(x)(dx^4)^2,
\label{general-foliation}
\ee
where $\Phi(x)$, $\beta_i(x)$, $\Psi(x)$ are functions at least $C^3$ 
in the whole manifold.
Clearly the above line element represents 
the most general four dimensional pseudo-Riemannian manifold
admitting a local foliation \footnote{By local we mean 
a foliation defined in a neighbourhood of a point which is 
covered by a coordinate chart.} 
by three dimensional conformally flat 
Riemannian hypersurfaces (given by the condition $x^4=const$). 
Now if we consider the integrable distribution associated to this 
foliation we may define an orthogonal projector 
$P_{ab}$ by means of (\ref{distribution}). 
To achieve this, we take the vector fields $\{\d/\d x^1,\d/\d x^2,\d/\d x^3\}$
which span the afore-mentioned distribution and construct an orthonormal 
set out of them (they are already orthogonal so it is enough with normalising).
Applying (\ref{distribution}) we easily get
$$
P_{11}=P_{22}=P_{33}=\Phi(x),\ P_{i4}=\beta_i(x),\ i=1,2,3,\ P_{44}=\sum_{i=1}^3\beta^2_i(x)/\Phi(x).
$$
Using this we can check condition (\ref{rank-3}) using the projector
$P_{ab}$ to calculate $L^0_{ab}$ and see what is obtained. 
The result is that the tensor $\bnb_{[a}L^0_{b]c}$ does not vanish in this 
case although a calculation  shows the important property
\be
P^r_{\ a}P^s_{\ b}P^q_{\ c}\bnb_{[r}L^0_{s]q}=0.
\label{changed-condition}
\ee
\begin{thm}
A necessary condition that a four dimensional 
pseudo-Rieman\-ni\-an manifold can be foliated locally by 
conformally flat Riemannian hypersurfaces is equation 
(\ref{changed-condition}).
\end{thm}  
 This result suggests 
that it may well be possible to generalize the conditions of theorem \ref{case-3} to pseudo-Riemannian manifolds of arbitrary 
dimension which are not conformally separable replacing these conditions 
by (\ref{changed-condition}). In fact 
this example can be generalized if we consider pseudo-Riemannian manifolds 
of higher dimension locally foliated by conformally flat hypersurfaces 
of dimension arbitrary and not necessarily Riemannian.
The condition which must be checked in this case is 
$T(P)^a_{\ bcd}=0$ being this condition found to be true in 
all the examples tried. Therefore it seems that a simple modification of 
the conditions 
of theorems \ref{bi-conformal-char-2} and 
\ref{case-3} holds even though the pseudo-Riemannian manifold is not 
conformally separable.  The true extent of this assertion is 
under current research. 
\end{exmp}

\section*{Acknowledgements}
We wish to thank Jos\'e M. M. Senovilla for a careful reading of 
the manuscript and his suggested improvements.
The valuable comments of an annonymous referee are also appreciated. 
Finally, financial support from projects\\  9/UPV00172.310-14456/2002 
and FIS2004-01626 is also gratefully acknowledged.

\end{document}